# Hadwiger's conjecture for 8-coloring graph

T.-Q. Wang[1] & X.-J. Wang[2]

[1]*Department of Mechanical and Manufacturing Engineering, Miami University, Oxford, Ohio, 04560, USA*

[2]*State Key Laboratory of Applied Optics, Changchun Institute of Optics & Fine Mechanics and Physics, Chinese Academy of Sciences, Changchun, Jilin, 130033, China*

**Hadwiger's Conjecture has been an open problem for over a half century[1,6], which says that there is at most a complete graph $K_t$ but no $K_{t+1}$ for every t-colorable graph. A few cases of Hadwiger's Conjecture, such as 1, 2, 3, 4, 5, 6-colorable graphs have been completely proved to convince all[1-5], but the proofs are tremendously difficult for over the 5-colorable graph[6,7]. Although the development of graph theory inspires scientists to understand graph coloring deeply, it is still an open problem for over 7-colorable graphs[6,7]. Therefore, we put forward a brand new chromatic graph configuration and show how to describe the graph coloring issues in chromatic space. Based on this idea, we define a chromatic plane and configure the chromatic coordinates in Euler space. Also, we find a method to prove Hadwiger's Conjecture for every 8-coloring graph feasible.**

Hadwiger's Conjecture, any graph for any positive integer t can be contracted at most to the complete graph with t-vertices if it can be t-colored, has been an open problem since 1943[1,4,6,7,17], which has been considered to be a huge strengthening of the four-color theorem[6,7].

For t =1, 2, 3, Hadwiger's Conjecture is true[1]. For t = 4, Hadwiger[1] and Dirac[3] proved it is workable. For t = 5, it becomes hugely difficult to attack. But Wagner proved that this case was equivalent to the Four Color Theorem in 1937[3]. Up to now, the four-color theorem: every planar graph is 4-colorable, has been proved, which proofs by use of computers were in 1976[8-10], 1997[11] and 2005[12], especially in developing to be simple and understand, but no hand proof, which is a 'humanoid proof'. Thus, for j = 5, it is true.

For t = 6, Seymour and Thomas proved it is true[4], which is one of the deepest in attacking this conjecture. For t = 7, Kawarabayashi and Toft proved that any 7-chromatic graph has $K_7$ or $K_{4,4}$ as a minor[5]. This may not be complete.

For t > 7, it is an open problem.

Here, every discussed graph is finite, undirected, simple, which without loops and multiple edges; $K_n$ is the complete graph with n vertices. Some graphic concepts are defined in the materials[13,14].

One of the milestones in the development of graph coloring is the concept: *minor*. The minor is a complete graph after edge contractions, in which the distinct but necessary colors are applied to color all adjacent vertices in distinct colors for a graph[15,16].

Based on the four-color theorem[8-10] and developing graph minor theorem from 1983 to 2010[15,16], we configure the new model for graph coloring as follows.

A chromatic plane (CHP) is a plane in 3-D Euler space. We stipulate that a chromatic subgraph (CHPsG) on a CHP is in at most four colors, which graph is the coloring subgraph of a graph on a CHP; The colors used in coloring CHPsG on each CHP are different from that on the other CHPs.

According to workable Hadwiger's Conjecture (for t < 6), the minor for every CHPsG is at most $K_4$ or every CHPsG on each CHP is four-colorable, in which the distinct colors used in a CHP is called CHP chromatic number that is at most four.

Besides, CHPsG may include disconnected subgraphs of a graph, but, for every CHP, it only allows putting on the disconnected graphs with the same minor.

A non-plane edge is an edge connecting vertices in two different CHPs, and a plane edge is an edge connecting vertices in the same CHP.

A chromatic-space (CH-3D), or chromatic coordinates, combines a series of parallel chromatic-planes in 3-D Euler space, in which a graph in every CHP is connected to all others in the other CHPs by non-plane edges.

The colors on each CHP, equaling chromatic number, is the third dimension in CH-3D.

A chromatic pallet paralleled to CHPs is a plane in which CHPsGs in several CHPs are projected or plotted together with the connection edges and vertices. The projecting operation keeps any vertex not covered by the others and connects all non-plane edges on the pallet to their belonging vertices.

To describe how to keep the minor of every CHPsG unchanged, we confine the generalized edge contractions for the graph.

*In a CH-3D, the edge contraction only contracts the plane edges by removing the edge, merging connected vertices into a new one for inheriting the other connecting edges connected to two merging vertices, including the non-plane edges.*

**Lemma 1.** *The minor of a graph keeps unchanged after edge contractions in a CH-3D.*

**Proof 1.** *If projecting all CHPsGs on every CHP into one chromatic pallet, the operation in Lemma 1 is one part of the generalized edge contractions. Thus Lemma 1 holds and the minor keeps unchanged.*

*But, if removing the non-plane edge, the minor is changed. One counter-example of removing non-plane edge operation is $K_5$. We put $K_5$ in two CHPs, minor $K_4$ in CHP(1), vertex $V_5$ in CHP(2). If removing a non-plane edge in red in Fig. 1, then, its minor becomes $K_4$.*

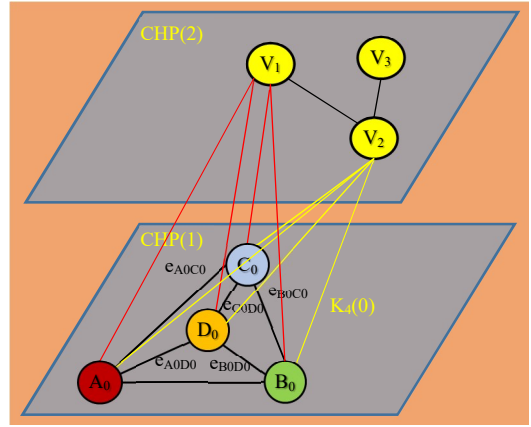

**Figure 1** Counter-example for removing non-plane edge in $K_5$. If any red line is removed, CHP chromatic number becomes 4+0, $V_1$ should be put on CHP(1). So, colors on CHP(2) are changed.

To recognize a vertex connected to the CHP whether or not on the CHP, we discuss the discriminant criterion in Lemma 2 and Corollaries, as follows.

**Lemma 2.** *If the number of different colors in which the vertices connected to this vertex are colored, is no less than the CHP chromatic number, the vertex is not on the CHP.*

**Proof 2.** *According to CHP definition, CHPsG is t-colorable (t<5). We color all vertices of CHPsG in t distinct colors. If the vertex V connects to the vertices in t different colors, then it should be in color t+1 to differ from its adjacent vertices. However, the CHP only allows t colors. Inversely, V is not in this CHP.*

For example, supposing vertex $V_1$ is connected to the four-coloring vertices on the CHP(1) in Fig. 1, which CHP(1) chromatic number is four, we can not color V in these four colors which vertices connected to it are colored in. We need to color it in the fifth color. Since the CHP(1) chromatic number is four, we are not allowed to put it on, then, we put this vertex on the other CHP such as CHP(2) in Fig. 1.

**Corollary 2.1.** *If a vertex is connected to vertices on the CHP, the connecting edges are*

*all-or-nothing either non-plane or plane edges of the CHP.*

**Proof 2.1.** *When a vertex is connected to the vertices on the CHP, there are two cases: if the vertex is on the CHP, the connected edges are plane ones on the CHP; if the vertex is not on the CHP, the connected edges are non-plane ones, which the number of non-plane edges connected to different coloring vertices is no less than the colors used on the CHP by using Lemma 2. Thus, Corollary 2.1 holds.*

**Corollary 2.2.** *If a vertex is connected to vertices on the CHP with non-plane edges, then it is completely connected to the minors of the CHP after edge contradictions.*

**Proof 2.2.** *By using Lemma 2 and Corollary 2.1, if the vertex V is connected to the vertices on the CHP by the non-plane edges, then these non-plane edges are connected to available vertices that are in at least including all distinct colors on the CHP.*

*By using Lemma 1, those non-plane edges are not annihilated but merged to the coloring vertices of the minor on the CHP, which minor is produces after edge contractions. Since the coloring minor of the CHP is colored in all available different colors due to workable Hadwiger's Conjecture (t<7). Thus, Corollary 2.2 holds.*

**Corollary 2.3.** *If the vertex is connected to the vertices on the same CHP, then it can be colored in at least one of the available colors of the CHP.*

**Proof 2.3.** *By using Lemma 2 and Corollary 2.2, if the vertex is the vertices on the same CHP, the number of its edges connected to different coloring vertices is less than the colors used on the CHP. Then it can be colored in the residual colors of the available colors of the CHP. Thus, Corollary 2.3 holds.*

**Corollary 2.4.** *The minors of every CHPsG are completely connected if without disconnected subgraphs on every CHPsG.*

**Proof 2.4.** *If without disconnected minors on every CHPsG, only one minor of every CHPsG is on each CHP. According to Corollary 2.2, each vertex of the minor on each CHP is completely connected to the other vertices of the minors on the other CHPs. Then, Corollary 2.4 holds.*

**Corollary 2.5.** *The graph of all minors on each CHP is a complete graph if without disconnected subgraphs on every CHPsG.*

**Proof 2.5.** *By using Corollary 2.4, without the connected edge, all minors are completely connected, so one of the vertices in all minors is all connected to the others. Then, Corollary 2.5 holds.*

**Corollary 2.6.** *To compare two coloring graphs, the larger the minor, the larger the chromatic number.*

**Proof 2.6.** *The size of the minor is the necessary colors needed for a coloring graph. The larger minor is used more colors than another.*

Based on defined edge contractions in the CH-3D in Lemma 1 and how to recognize each vertex whether or not on the CHP in Lemma 2, we find a way, chromatic-filling CHPs for a graph, to assign a graph to CHPs in CH-3D to keep without disconnected graph in each CHPs and separating edge contractions from non-plane edges.

According to the adjacent list of an 8-coloring graph G, (1), do an initiation and let subG(1) = G and CHPsG(1) is empty.

(2), find the maximal minor[4] $K_4$ in graph subG(1), we at will select a subgraph which minor is as same as $K_4$, and put this subgraph into CHPsG(1), and color the vertices of CHPsG(1) in 4 distinct colors, and form the 1st CHP with the chromatic number 4 as CHP(1);

(3), enumerate all vertices in the residual graph as ResidG(1) = subG(1)-CHPsG(1);

(4), according to Lemma 2, if the vertex is not connected to CHP(1) and is not on CHP(1) by using the adjacent list of the graph, keep it in the residual graph ResidG(1).

(5), if the vertex is on CHP(1), add it in CHPsG(1), and color it in the residual color from one of the four colors of CHP(1); rebuild the dynamical ResidG(1) = subG(1)-CHPsG(1); re-enumerate ResidG(1) by using the adjacent list; return to step (3).

(6), at the end of the enumeration of the ResidG(1), if ResidG(1) is empty, then the graph is all assigned to the CHPs in a CH-3D.

(7), otherwise, if disconnected graphs appear in ResidG(1), we select the disconnected graph with the maximal minor as the new iterative subG (2) according to Corollary 2.6; return to step (2).

(8), after assigning all vertices for the at most coloring graph in step (6), we can put on non-plane and plane edges depending on the adjacent list of the graph.

**Lemma 3.** *Every finite 8-coloring graph with the maximal minor can be assigned in CH-3D after color-filling operations.*

**Proof 3.** *By using chromatic-filling CHPs, every vertex is checked out its relation to its connected CHPs according to the adjacent list of a graph.*

*Moreover, the disconnected graphs might be created in each residual graph after each iteration. To keep the iterative subgraph is connected, the criterion is to remove the disconnected graphs which minors are less than the maximal minor among disconnected graphs. Then, we substitute the disconnected graphs with the maximal minor among disconnected graphs for the residual graph in the iteration to ensure we can assign the graph with the maximal minor to CHPs in CH-3D. Since the graph is finite, the iteration is limited and implementable.*

*Therefore, to describe a coloring graph in CH-3D, we can assign its vertices in the graph with the maximal minor to each CHPsG on each CHP, including connected non-plane and plane edges.*

**Corollary 3.1.** *In chromatic-filling CHPs for a graph, if a cut edge or a cut vertex is assigned on the CHP, and the minors of the disconnected graphs linked by a cut edge or a cut vertex in the residual graph for the next iteration are more than the chromatic number of the CHP, the residuals graph has at least one disconnected graph at the end of the enumeration.*

**Proof 3.1.** *Supposing there is one cut edge or one cut vertex that appears in subG(i) for assigning subG(i) on the i-th CHP(i), the cut edge or vertex divides subG(i) to at least the left and right hand, the left hand will be residual subgraph L and so as subgraph R for the right hand. If the minors of subgraph L and R are more than the chromatic number of the CHP(i) and the one cut edge or one cut vertex is assigned on the CHP(i), at the end of the enumeration, subgraph L is disconnected to subgraph R. So, the residuals graph ResidG(i) has at least one disconnected graph.*

**Corollary 3.2.** *The larger minor between the disconnected L and R decides the necessary colors for the coloring graph at most.*

**Proof 3.2.** *We can color the disconnected graphs independently. The graph with large minor needs more colors than that with less minor according to Corollary 2.6.*

**Corollary 3.3.** *All minors for each CHPs are completely connected and form a complete graph after chromatic-filling operations.*

**Proof 3.3.** *By using Lemma 2, every vertex is connected to the others with either non-plane edges or plane edges. After chromatic-filling operations, even if a vertex is on a CHP, it is connected to vertices of CHPsG on the same CHP with at least one edge.*

*After the chromatic-filling CHPs, all minors for each CHPsG are completely connected by using Corollary from 2.1 to 2.6, since we replace the possible disconnected graph for one with the maximal minor. So, all minors form a complete graph.*

**Hadwiger's Cojecture[17] at t=8.** *for every 8-*

*coloring graph, it has at most $K_8$.*

**Proof.** Basic configuration for every 8-colorable graph, it has double $K_4$ as minors after chromatic-filling CHPs in a CH-3D. Since every 8-colorable graph is colored 8 colors at most, we find a subgraph with the minor[4] as $K_4$ assigning it in the subG(1) for the iterations in the chromatic-filling operations.

By using chromatic-filling CHPs to assign the graph into CHPs, CHPsG(1) is obtained after the 1st iteration; Then, CHPsG(2) after the 2nd iteration; Until the residual set is empty, the iteration is stopped, according to Lemma 3.

Especially, vertex $V_1$ is the one in the 8-colorable graph, there are at least four edges for $V_1$ connecting to four distinct coloring vertice on CHP(1). So, $V_1$ is not on CHP(1) in the 1st iteration. Since CHPsG(1) is removed from the graph before the 2nd iteration, $V_1$ is connected to only at least three distinct coloring vertices on CHP(2). So, $V_1$ is on CHP(2).

According to feasible Hadwiger's Conjecture (t<7), we can obtain the minor for each CHP by using Lemma 1, shown in Fig. 2.

According to Corollary 3.3, we obtain double $K_4$ that is completely connected. Thus, we obtain $K_8$ at last. Hadwiger's Conjecture at t=8 holds.

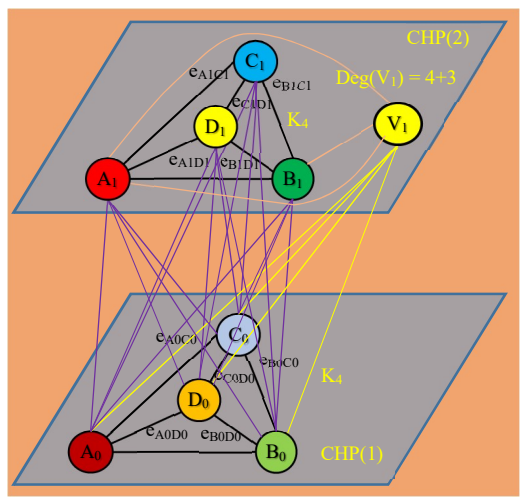

**Figure 2** Schematic graph for Hadwiger's Conjecture: t = 8 for the graph having double $K_4$ as minors.

On another hand, if we know the graph with minor $K_9$, then we can assign the graph to two CHPs, which the first CHP is with chromatic number 4. The second CHP is with chromatic number at most 4. So, the graph is 8-colorable at most according to Corollary 3.3.

As above, we proved that Hadwiger's Conjecture for an 8-colorable graph holds by using the configuration of chromatic coordinates.

We are enlightened by the former's proof to understand the coloring graph deeply. The defined chromatic plane is efficient to solve the contradiction between graph coloring and plane graph, and chromatic coordinates can visually split interwinding coloring graphs.

We utilize the color-filling operation assigning every coloring graph in chromatic coordinates and configure the minors in each chromatic plane are completely connected, which is conclusive for the proof on Hadwiger's Conjecture.

We hope that concepts of chromatic coordinates and operations can play more and more important roles in applied sciences.

**Author's contribution** All authors equally contributed to design, proving, and writing the manuscript.

**Acknowledgements** T.W. acknowledges support from the fellowship of MME, Miami University, Oxford. X.W. thanks the Foundation of SKLAO, Innovation Foundation of CAS, and the NSF of China.

**Competing interests statement** The authors declare that they have no competing financial interests.

**Correspondence** and requests for materials should be addressed to T.W. (wangt5@miamioh.edu) or X.W. (xjwang@ciomp.ac.cn)